\documentclass{article}

\usepackage[T1]{fontenc}
\usepackage{amsmath}
\usepackage{amsfonts}
\usepackage{verbatim}
\usepackage{bbm}
\usepackage{amsthm, amssymb}

\newtheorem{theorem}{Theorem}

\newtheorem{proposition}[theorem]{Proposition}

\theoremstyle{definition}

\theoremstyle{remark}
\newtheorem{remark}[theorem]{Remark}

\theoremstyle{remark}
\newtheorem{example}[theorem]{Example}

 \numberwithin{equation}{section}
\numberwithin{theorem}{section}

\date{}

\begin{document}

\title{Asymptotic Normality of Random Sums of
m-dependent Random Variables}

\author{ \"{U}mit I\c{s}lak}

\maketitle

\begin{abstract}
We prove a central limit theorem for random sums of the form
$\sum_{i=1}^{N_n} X_i$,  where $\{X_i\}_{i \geq 1}$ is a stationary
$m-$dependent process and $N_n$ is a random index independent of
$\{X_i\}_{i\geq 1}$. Our proof is a generalization of Chen and
Shao's result
 for i.i.d. case and consequently
we recover their result. Also a variation of a recent result of
Shang on $m-$dependent sequences is obtained as a corollary.
Examples on moving averages and descent processes are provided, and
possible applications on non-parametric statistics are discussed.
\end{abstract}

\section{Introduction }

In the following, we analyze the asymptotic behavior of random sums
of the form $\sum_{i=1}^{N_n} X_i$ as $n \rightarrow \infty$, where
$X_i'$s are non-negative random variables that are stationary and
$m-$dependent, and $N_n$ is a non-negative integer valued random
variable independent of $X_i's$. Limiting distributions of random
sums of independent and identically distributed (i.i.d.) random sums
are well studied. See \cite{chenshaopaper}, \cite{klaver},
\cite{robbins} and the references therein. Asymptotic normality of
deterministic sums of $m-$dependent random variables are also well
known. See, for example, \cite{bergstrom}, \cite{hoeffding} and
\cite{orey}.  To the best of author's knowledge, previous work on
the case of random sums of the form $\sum_{i=1}^{N_n} X_i$  where
$X_i$'s are dependent are limited to \cite{shang} where he works on
$m-$dependent random variables and \cite{barbour} where they
investigate random variables that appear as a result of integrating
a random field with respect to  point processes.  Our results here
will be in the lines of \cite{chenshaopaper} generalizing their
result to the stationary $m-$dependent case. Throughout the way, we
will also improve the results given in \cite{shang}.

Let's now recall  stationary and $m-$dependent processes. Let
$\{X_i\}_{i\geq 1}$ be a stochastic process and let
$F_X(X_{i_1+m},...,X_{i_k+m})$ be the cumulative distribution
function of the joint distribution of $\{X_i\}_{i\geq 1}$ at times
$i_1+m,...,i_k+m$. Then $\{X_i\}_{i\geq 1}$ is said to be
\emph{stationary } if, for all $k$, for all $m$ and for all
$i_1,...,i_k$
$$F_X(X_{i_1+m},...,X_{i_k+m})=F_X(X_{i_1},...,X_{i_k})$$
holds. For more on stationary processes, see \cite{shiryaev}. If we
define the distance between two subsets of $A$ and $B$ of
$\mathbb{N}$ by $$\rho(A,B):= \inf \{|i-j|: i \in A, j \in B\},$$
then the sequence $\{X_i\}_{i\geq 1}$ is said to be
\emph{$m-$dependent} if $\{X_i, i \in A\}$ and $\{X_j, j \in B\}$
are independent whenever $\rho(A,B)
> m$ for $A, B \subset \mathbb{N}.$

An example of a stationary $m-$dependent process can be given by the
moving averages process. Assume that $\{T_i\}_{i \geq 1}$ is a
sequence of i.i.d. random variables with finite mean $\mu$ and
finite variance $\sigma^2.$ Letting $X_i =(T_i + T_{i+1})/2$,
$\{X_i\}_{i \geq 1}$ is a stationary 1-dependent process with
$\mathbb{E}[X_i] = \mu,$ $Var(X_i ) = \sigma^2 /2 $ and
$Cov(X_1,X_2) = \sigma^2 /4.$

This paper is organized as follows: In the next section, we state
our main results and compare them with previous approaches. In the
third section, we give examples on moving averages and descent
processes relating it to possible nonparametric tests where the
number of observations is itself random. Proofs of the main results
are given in Section 4 and we conclude the paper with a discussion
of future directions.

\section{Main Results}

We start with two propositions. Proofs of these are standard and are
given at the end of Section \ref{Proofs}.

\begin{proposition}\label{var1}
Let $\{X_i\}_{i \geq 1}$ be a stationary $m-$dependent process with
$\mu:=\mathbb{E}[X_i]$, $\sigma^2 := Var(X_i)< \infty$,
$a_{j}:=Cov(X_1,X_{1+j})$. Then for any $N \geq 1$, we have
$$Var\left(\sum_{i=1}^N X_i \right) = N (\sigma^2 + 2 \sum_{j=1}^m a_j \Gamma_{N,j} )-2\sum_{j=1}^m j a_j \Gamma_{N,j}
$$ where $\Gamma_{N,j} = \mathbbm{1}(N \geq j+1).$
\end{proposition}

\begin{proposition}\label{varrandom}
Let $\{X_i\}_{i\geq 1}$ be as in Proposition \ref{var1}. Let $Y_i$'s
be i.i.d. non-negative integer valued random variables with $\nu:=
\mathbb{E}[Y_i]$, $\tau^2:= Var(Y_i)<\infty$ and assume that $X_i$'s
and $Y_i$'s are independent. Define $N_n = \sum_{i=1}^n Y_i.$ Then
we have
$$Var\left(\sum_{i=1}^{N_n} X_i \right) = n( \nu \sigma^2+ 2  \nu \sum_{j=1}^m a_j+ \mu^2 \tau^2) +
\alpha(m)$$ where $$\alpha(m) = \sum_{k=0}^m (2k \sum_{j=1}^m
a_j(\Gamma_{k,j}-1)-2 \sum_{j=1}^m j a_j (\Gamma_{k,j}-1)) -2
\sum_{j=1}^m j a_j$$ and $\Gamma_{k,j} = \mathbbm{1}(k \geq j+1)$.
In particular, $\frac{\alpha(m)}{n} \longrightarrow 0$ as $n
\rightarrow \infty.$ When $X_i's $ are also independent (i.e.,
$m=0$), this reduces to $$Var\left(\sum_{i=1}^{N_n} X_i \right) = n
(\nu \sigma^2 + \mu^2 \tau^2).$$
\end{proposition}

In the following, we will be using $\longrightarrow_d$ for
convergence in distribution and $=_d$ for equality in distribution.
Also $N(0,1)$  and $\Phi$ will denote a standard normal random
variable and its cumulative distribution function, respectively. Now
we are ready to present our main result.

\begin{theorem}\label{mainresult} Let $\{X_i\}_{i \geq 1}$ be a non-negative stationary
$m-$dependent process with $\mu:=\mathbb{E}[X_1] >0$, $\sigma^2 :=
Var(X_1) > 0$, $a_{j}:=Cov(X_1,X_{1+j})$, $\sigma^2+2 \sum_{j=1}^m
a_j
>0$ and  $\mathbb{E}|X_1|^3 < \infty.$ Let $Y_i$'s be
i.i.d. non-negative integer valued random variables with $\nu:=
\mathbb{E}[Y_1]>0$, $\tau^2:= Var(Y_1)>0$, $\mathbb{E}|Y_1|^3 <
\infty$ and suppose that $X_i$'s and $Y_i$'s are independent. Define
$N_n = \sum_{i=1}^n Y_i.$ Then
\begin{equation}\label{maintheorem}
\frac{\sum_{i=1}^{N_n} X_i - n\mu \nu}{\sqrt{n (\nu \sigma^2 + 2 \nu
\sum_{j=1}^n a_j + \tau^2 \mu^2)}} \longrightarrow_d N(0,1)
\end{equation}
as $n \rightarrow \infty$.
\end{theorem}

Note that assumptions on $Y_i$'s hold, for example, when $Y_i$'s are
 non-degenerate i.i.d. Bernoulli random variables. This is one of the most natural
 cases as in that case we may consider $\sum_{i=1}^{N_n} X_i$ as the sum of
outcomes of a series of experiments, where each observation is
blocked with a fixed probability independent of others. The main
assumption on $X_i'$s (others are non-degeneracy conditions) is a
third moment condition.

 Since our proof is a direct
generalization of Chen and Shao's result on i.i.d. case (which is
the case with $m=0$), we recover their result from
\cite{chenshaopaper}.

\begin{theorem}\label{chencase}
Let $\{X_i\}_{i \geq 1}$ be i.i.d. random variables with
$\mu:=\mathbb{E}[X_1]
>0$, $\sigma^2 := Var(X_1)>0$,
 and assume that $\mathbb{E}|X_1|^3 < \infty.$ Let $Y_i$'s be
i.i.d. non-negative integer valued random variables with $\nu:=
\mathbb{E}[Y_1]>0$, $\tau^2:= Var(Y_1)>0$, $\mathbb{E}|Y_1|^3 <
\infty$ and assume that $X_i$'s and $Y_i$'s are independent. Define
$N_n = \sum_{i=1}^n Y_i.$ Then for any $n\geq 1$, we have
\begin{equation}\label{chensbound}
 \sup_{z \in \mathbb{R}} \left| \mathbb{P} \left( \frac{\sum_{i=1}^{N_n} X_i - n\mu \nu}{\sqrt{n (\nu \sigma^2  + \tau^2 \mu^2)}} \leq z\right) - \Phi(z)\right|
\leq C n^{-1/2} \left(\frac{\tau^2}{\nu^2} +
\frac{\mathbb{E}[Y_1^3]}{\tau^3} +
\frac{\mathbb{E}|X_1|^3}{\nu^{1/2} \sigma^3} + \frac{\sigma}{\mu
\sqrt{\nu}} \right)
\end{equation}

\noindent where $C$ is a constant independent of $n$.
\end{theorem}

We will explain how the proof of Theorem \ref{mainresult} also
reveals Theorem \ref{chencase} in Section \ref{Proofs}. We note that
in the original statement of Chen and Shao's result, $\mu$ is
allowed to be 0. We excluded this in our statement as the upper
bound in (\ref{chensbound}) is $\infty$ when $\mu =0$.

Our final result will be a variation of  the main theorem given in
\cite{shang} about the asymptotics of random sums of $m-$dependent
random variables. Namely, we have

\begin{theorem}\label{shang}
Under the assumptions of Theorem \ref{mainresult},
\begin{equation}\label{shangtheorem}
\frac{\sum_{i=1}^{N_n} X_i - N_n \mu}{\sqrt{N_n} \left(\sigma + 2
\sum_{j=1}^m a_j \right)} \longrightarrow_d N(0,1)
\end{equation}
 as $n \rightarrow
\infty.$
\end{theorem}

\begin{remark}
Indeed, as can be seen from the proof of Theorem \ref{mainresult},
one can obtain convergence rates when the scaling is perturbed a
little bit. More precisely, we have $$\sup_{z \in \mathbb{R}} \left|
\mathbb{P}\left( \frac{\sum_{i=1}^{N_n} X_i - N_n \mu}{\sqrt{N_n
\sigma'} } \leq z\right) - \Phi(z) \right| \leq \frac{C}{\sqrt{n}}$$
for a universal constant $C$ and for every $n \geq 1$ where
$(\sigma')^2 = \sigma^2 + 2 \sum_{j=1}^m a_j \Gamma_{N_n,j} -
\frac{2}{N_n} \sum_{j=1}^m j a_j \Gamma_{N_n,j}$ and $\Gamma_{N_n,j}
=\mathbbm{1} (N_n \geq j+1)$.
\end{remark}

\section{Examples}

\begin{example} (Moving averages)
Assume that $\{T_i\}_{i \geq 1}$ is a sequence of i.i.d. random
variables with finite mean $\mu$ and finite variance $\sigma^2.$
Letting $$X_i =\frac{T_i + T_{i+1}}{2}, \; \; i\geq 1,$$ $\{X_i\}_{i
\geq 1}$ is a stationary 1-dependent process with $\mathbb{E}[X_i] =
\mu,$ $Var(X_i ) = \sigma^2 /2 $ and $Cov(X_1,X_2) = \sigma^2 /4.$
When $\mu >0$  and $\sigma^2 >0$, we can apply Theorem
\ref{mainresult} as long as the assumptions on $N_n$ are satisfied
(As noted above, they will be satisfied, for example, when $Y_i$'s
are independent Bernoulli random variables with success probability
$p \in (0,1)$). This discussion can be generalized to $m-$moving
averages defined as
$$Y_i = \frac{T_i+T_{i+1}+...+T_{i+m-1}}{m}, \quad i \geq 1, \quad m \in \mathbb{N}$$
in a straightforward way.
\end{example}

\begin{example} (Descent processes)
A sequence of real numbers $(t_i)_{i=1}^n$ is said to have a
\emph{descent} at position $1 \leq k \leq n-1$ if $t_k > t_{k+1}$.
 Here we are interested in the descent process of a sequence
of random variables.  Statistics related to descents
 are often used in nonparametric statistics to test independence or
 correlation (For example, one uses the number of inversions in Kendall's tau statistic).
  See \cite{ferguson} for a brief introduction for this
 connection. Also see \cite{diaconis}  to learn more about why
 these processes are important.

Now let $T_i$'s be i.i.d. random variables with distribution $F$,
and $X_i:=\mathbbm{1}(T_i
> T_{i+1}).$ Also let $Y_i'$s be i.i.d. Bernoulli random variables with parameter $p \in (0,1)$ and set $N_n = \sum_{i=1}^{N_n} Y_i$. Defining $$W_n = \sum_{j=1}^{N_n-1} X_i,$$
 $W_n$ is
the number of descents in the random length sequence
$(T_1,T_2,...,T_{N_n}).$

Here $\{X_i\}_{i \geq 1}$  is a stationary 1-dependent process and
it is easy to check that $\mu = 1/2, \; \sigma^2 = 1/4 \quad
\text{and} \quad \sigma^2+2a_1 = 1/12.$ So assumptions of  Theorem
\ref{mainresult} are satisfied and we obtain the asymptotic
normality of $W_n$.
\end{example}

\begin{example} (Non-parametric statistics ) In this example, we discuss a possible
application of Theorem \ref{mainresult} in non-parametric
statistics. Let $T_1,...,T_n$ be the random outcomes of an
experiment and assume that the probability of observing any of these
is $p \in (0,1)$ independent of each other. Let $N_n$ be the number
of actually observed outcomes and $O_1,...,O_{N_n}$ be the
corresponding sequence of observations.

Suppose we want to test $$H_0: T_1,...,T_n \quad \text{are
uncorrelated and } \quad p=p_0. $$ Then one can use the test
statistic
$$W_n = \sum_{i=1}^{N_n-1} \mathbbm{1} (O_i > O_{i+1})$$ and Theorem \ref{mainresult} to understand
the asymptotic distribution of $W_n$ under the null hypothesis. A
very large or a very small value for this statistics will provide
information about the dependence structure of $T_i'$s.

Extensions of this observation to more general tests will be
followed in a subsequent work.
\end{example}

\section{Proofs}\label{Proofs}

We start by recalling two results that will be useful in the proof
of the main theorem. First of these is a central limit theorem for
$m-$dependent random variables established in \cite{chenshaolocal}.

\begin{theorem}\label{localdepresult} \cite{chenshaolocal} If $\{X_i\}_{i \geq 1}$ is  a sequence of zero mean
$m-$dependent random variables and $W=\sum_{i=1}^{n} X_i$, then for
all $p \in (2,3]$, $$\sup_{z \in \mathbb{R}} |\mathbb{P}(W \leq z) -
\Phi(z)| \leq 75 (10m+1)^{p-1} \sum_{i=1}^{n} \mathbb{E}|X_i|^p.$$

\end{theorem}

The second result we will need is  the following theorem of Chen and
Shao (\cite{chenshaopaper}) on the normal approximation of random
variables. We note that this theorem is part of what is known as the
concentration inequality approach in Stein method literature. See
the cited paper or \cite{goldsteinbook} for more on this.

\begin{theorem}\label{concentration} \cite{chenshaopaper} Let $\xi_1,...,\xi_n$ be independent
mean zero random variables for $i=1,...,n$ with $\sum_{i=1}^n
Var(\xi_i)=1$. Let  $W= \sum_{i=1}^n \xi_i$, $T= W + \Delta$, and
also for each $i = 1,...,n,$ let $\Delta_i$ be a random variable
such that $\xi_i$ and $(W-\xi_i,\Delta_i)$ are independent. Then we
have
\begin{equation}\label{chenshaobound}
\sup_{z \in \mathbb{R}} |\mathbb{P}(W \leq z) - \Phi(z)| \leq
6.1(\beta_2 + \beta_3) + \mathbb{E}|W \Delta| + \sum_{i=1}^n
\mathbb{E} |\xi_i (\Delta - \Delta_i)|
\end{equation}

\noindent where
$$\beta_2 = \sum_{i=1}^n \mathbb{E} [\xi_i^2 \mathbbm{1}(|\xi_i|>1)]
\qquad \text{and} \qquad \beta_3 = \sum_{i=1}^n \mathbb{E}
[|\xi_i|^3 \mathbbm{1}(|\xi_i|\leq 1)].$$
\end{theorem}

Before moving on to the proof of Theorem \ref{mainresult},  we
finally  recall Prokhorov and Kolmogorov distances between
probability measures. Let $\mathcal{P}(\mathbb{R})$ be the
collection of all probability measures on
$(\mathbb{R},\mathfrak{B}(\mathbb{R}))$ where
$\mathfrak{B}(\mathbb{R})$ is the Borel sigma algebra on
$\mathbb{R}$. For a subset $A \subset \mathbb{R},$  define the
$\epsilon-$neighborhood of $A$ by
$$A^{\epsilon} := \{p \in \mathbb{R}: \exists q \in A, d(p,q) < \epsilon\}= \bigcup_{p \in A}
B_{\epsilon}(p)$$ where $B_{\epsilon}(p)$ is the open ball of radius
$\epsilon$ centered at $p.$ Then the Prokhorov metric $d_p:
\mathcal{P}(\mathbb{R})^2\longrightarrow [0,\infty)$ is defined by
setting the distance between two probability measures $\mu$ and
$\nu$ to be

\begin{equation}\label{prokhorov}
    d_p(\mu,\nu):= \inf \{\epsilon >0 : \mu(A) \leq \nu(A^{\epsilon})+\epsilon \quad \text{and} \quad \nu(A) \leq \mu(A^{\epsilon})+\epsilon, \; \, \forall A \in \mathfrak{B}(\mathbb{R})\}.
\end{equation}

The Kolmogorov distance $d_K$ between two probability measures $\mu$
and $\nu $ is defined to be $$d_K(\mu, \nu) = \sup_{z \in
\mathbb{R}}|\mu((-\infty,z]) - \nu((-\infty,z]) |.$$

The following two facts will be useful: (1) Convergence of measures
in Prokhorov metric is equivalent to the weak convergence of
measures. (2) Convergence in Kolmogorov distance implies convergence
in distribution, but the converse is not true. See, for example,
\cite{shiryaev} for these standard results.

Now we are ready to prove Theorem \ref{mainresult}. We will follow
the notations of \cite{chenshaopaper} as much as possible.

\vspace{0.15in}

\textbf{Proof of Theorem \ref{mainresult} : }Let $Z_1, Z_2$ and
$Z_3$ be independent standard normal random variables which are also
independent of $X_i$'s and $Y_i$'s. Put
$$b= \sqrt{\nu \sigma^2 + 2 \nu \sum_{j=1}^m a_j + \tau^2 \mu^2}.$$
Define $$T_n = \frac{\sum_{i=1}^{N_n} X_i - n \mu \nu}{\sqrt{n} b}
\quad \text{and } \quad H_n = \frac{\sum_{i=1}^{N_n} X_i - N_n
\mu}{\sqrt{N_n} \sigma'}$$ where

\begin{equation}\label{sigmaprime}
    (\sigma')^2 = \sigma^2 + 2 \sum_{j=1}^m a_j \Gamma_{N_n,j} -
\frac{2}{N_n} \sum_{j=1}^m j a_j \Gamma_{N_n,j}
\end{equation}

\noindent  with $\Gamma_{N_n,j}:= \mathbbm{1}(N_n \geq j+1).$ Also
write
$$T_n = \frac{\sqrt{N_n} \sigma'}{\sqrt{n}b}H_n + \frac{(N_n - n\nu)
\mu}{\sqrt{n}b}$$ and
$$T_n(Z_1) =\frac{\sqrt{N_n}
\sigma'}{\sqrt{n}b}Z_1 + \frac{(N_n - n\nu) \mu}{\sqrt{n}b}.$$ For
$n$ large enough, we have $m < n \nu /2$. For such $n$, we have

\begin{eqnarray}\label{kolmbound}
  d_K(T_n,T_n(Z_1)) &=& d_K(H_n,Z_1) \notag \\
                   &\leq& \mathbb{P}(|N_n-n\nu|> n\nu/2) + \sup_{z \in
                   \mathbb{R}} \mathbb{E}\big[\mathbb{E}[|\mathbbm{1}(H_n \leq z)- \mathbbm{1}(Z_1 \leq z)| \mathbbm{1}(|N_n-n\nu| \leq
                   n\nu/2) \big|N_n]\big]  \notag\\
                   &\leq& \frac{4 \tau^2}{n \nu^2} + \mathbb{E}
                   \left[\mathbb{E}\left[\frac{C N_n \mathbb{E}|X_1|^3 \mathbbm{1}(|N_n-n\nu|\leq n\nu /2)}{N_n^{3/2} \left(\sigma^2 +2 \sum_{j=1}^m a_j - \frac{2}{N_n} \sum_{j=1}^m j a_j \right)^{3/2}}  \big| N_n\right]
                   \right]
\end{eqnarray}

\noindent where for (\ref{kolmbound}) we used Chebyshev's inequality
for the first estimate and Theorem \ref{localdepresult} with $p=3$
for the second estimate. Here the condition that $m< n\nu/2$
simplifies $(\sigma')^2$ as defined in (\ref{sigmaprime}) to
$(\sigma')^2 = \left(\sigma^2 + 2 \sum_{j=1}^m a_j - \frac{2}{N_n}
\sum_{j=1}^m j a_j \right)$ when $|N_n -n \nu| \leq n\nu /2$. Also
note that throughout this proof, $C$ will be a positive constant
with not necessarily the same value in different lines.  Now if
$\sum_{j=1}^m j a_j<0$, then the bound in (\ref{kolmbound}) yields

\begin{eqnarray*}
   d_K(T_n,T_n(Z_1)) \leq \frac{4 \tau^2}{n\nu^2}+ \frac{C\mathbb{E}|X_1|^3}{\sqrt{n \nu /2 } \left(\sigma^2+ 2 \sum_{j=1}^m a_j \right)^{3/2}}  \longrightarrow 0
\end{eqnarray*}

\noindent as $n \rightarrow \infty.$  Else if $\sum_{j=1}^m j a_j
\geq 0$, we observe that for large enough $n$, we have
$\sigma^2+2\sum_{j=1}^m a_j - \frac{4}{n \nu} \sum_{j=1}^m j a_j
>0$ by our assumption that $\sigma^2 + 2\sum_{j=1}^m a_j >0$. For such
$n$, using the bound in (\ref{kolmbound}) we obtain
\begin{eqnarray}\label{bound1}
  d_K(T_n,T_n(Z_1)) &\leq& \frac{4 \tau^2}{n\nu^2}+ \frac{C\mathbb{E}|X_1|^3}{\sqrt{n \nu /2 } \left(\sigma^2+ 2 \sum_{j=1}^m a_j - \frac{4}{n \nu} \sum_{j=1}^m j a_j\right)^{3/2}}
\end{eqnarray}

\noindent  and this yields   $d_K(T_n,T_n(Z_1)) \longrightarrow 0$
as $n \rightarrow \infty$  when $\sum_{j=1}^m j a_j \geq 0$.

Hence we conclude that $ d_K(T_n,T_n(Z_1)) \longrightarrow 0$ as $n
\rightarrow \infty$ as long as $\nu > 0$ and $\sigma^2 + 2
\sum_{j=1}^m a_j
>0.$ This in particular implies
\begin{equation}\label{tn1}
    d_p(T_n,T_n(Z_1)) \longrightarrow 0
\end{equation}
 as $n \rightarrow \infty$ where
$d_P$ is the Prokhorov distance as defined in (\ref{prokhorov}).

\vspace{0.2in}

Next let $(\sigma'')^2 = (\sigma')^2 + 2 \sum_{j=1}^m a_j (1 -
\Gamma_{N_n,j})+ \frac{2}{N_n} \sum_{j=1}^m j a_j \Gamma_{N_n,j}$ so
that $$(\sigma'')^2 = \sigma^2 + 2 \sum_{j=1}^m a_j.$$

Note that $\sigma''$ is not random and introduce
$$T_n'(Z_1)=\frac{\sqrt{N_n} \sigma''}{\sqrt{n}b}Z_1 + \frac{(N_n -
n\nu) \mu}{\sqrt{n}b}$$ and
$$T_n(Z_1,Z_2):=\frac{\tau \mu}{b} \left(Z_2+\frac{\sigma''
\sqrt{\nu}}{\tau \mu} Z_1
 \right).$$

One can easily check that $T_n(Z_1,Z_2)$ is a standard normal random
variable since $Z_1$ and $Z_2$ are assumed to be independent. So if
we can show that
 $d_p(T_n(Z_1),T_n'(Z_1)) \longrightarrow 0$ and $d_p(T_n'(Z_1),T_n(Z_1,Z_2)) \longrightarrow 0$ as $n
\rightarrow \infty$, then the result will follow from an application
of triangle inequality. We start by showing that $d_p(T_n'(Z_1),
T_n(Z_1,Z_2)) \rightarrow 0$.
 For this purpose, we will use Chen-Shao's concentration inequality
approach to get bounds in the Kolmogorov distance and to recover
Chen and Shao's result on i.i.d. case (If we just wanted to show
$d_P(T_n'(Z_1), T_n(Z_1,Z_2)) \rightarrow 0$, then this could be
done in a much easier way. See Remark \ref{prokhorovremark}). The
following argument is in a sense rewriting the corresponding proof
in \cite{chenshaopaper} with slight changes since the concentration
approach is used on $N_n$ which is in both problems a sum of
independent random variables. For the sake of completeness, we
include all details.

Define the truncation $\overline{x}$ of $x \in \mathbb{R}$ by

$$
\overline{x} = \begin{cases}
  n\nu/2 & \text{if $x < n  \nu /2$} \\
  x & \text{if $n\nu/2 \leq
x \leq 3n \nu /2$} \\
  3n \nu/2 & \text{if $x > 3 n \nu/2$}
\end{cases}$$

\noindent and let
$$\overline{T_n'}=\frac{\sqrt{\overline{N_n}}\sigma''}{\sqrt{n}b} Z_1+ \frac{(N_n-n\nu) \mu}{\sqrt{n} b} = \frac{\tau \mu}{b} \left(W + \Delta + \frac{\sigma'' \sqrt{\nu}}{\tau \mu}Z_1
\right)$$ where $$W= \frac{N_n-n\nu}{\sqrt{n}\tau} \quad \text{and}
\quad \Delta= \frac{(\sqrt{\overline{N_n}}-\sqrt{n \nu}) \sigma''
Z_1}{\sqrt{n} \tau \mu}.$$ Since $Y_i$ is independent of $N_n-Y_i$
for all $i=1,...,n$, we can apply Theorem \ref{concentration} to $W
+ \Delta $ setting
$$\Delta_i= \frac{\sqrt{\overline{N_n-Y_i+\nu}}-\sqrt{n \nu}\sigma''
Z_1}{\sqrt{n} \tau \mu}, \quad i=1,...,n.$$

(So $\xi_i =\frac{Y_i-\nu}{\sqrt{n} \tau}$ in Theorem
\ref{concentration}.) For the first term of the upper bound given in
(\ref{chenshaobound}), we have
\begin{equation}\label{bounds1}
6.1(\beta_2+\beta_3) \leq 6.1(2n) \mathbb{E}
\left|\frac{Y_1}{\sqrt{n} \tau} \right|^3 \leq \frac{Cn
\mathbb{E}|Y_1|^3}{(n\tau^2)^{3/2}}=\frac{C
\mathbb{E}|Y_1|^3}{\tau^3 \sqrt{n}}.\end{equation}

For the second term in (\ref{chenshaobound}), we have

\begin{eqnarray*}
  \mathbb{E}|W \Delta| &=&  \mathbb{E}|Z_1| \mathbb{E}\left[\frac{\sigma''}{\sqrt{n} \tau \mu} \mathbb{E}|W(\sqrt{\overline{N_n}}-\sqrt{n \nu})|\right] = \frac{\mathbb{E}|Z_1|\sigma''}{\sqrt{n} \tau \mu}
    \mathbb{E} \left|W \frac{\overline{N_n}-n\nu}{\sqrt{\overline{N_n}}+ \sqrt{n \nu}} \right|
\end{eqnarray*}
where we used the identity $\sqrt{x}-\sqrt{y} =
\frac{x-y}{\sqrt{x}+\sqrt{y}}$ in the second equality. So by an
application of Cauchy-Schwarz inequality, we obtain

\begin{eqnarray}\label{bounds2}
     \mathbb{E}|W \Delta| \leq  \frac{C\sigma''}{\sqrt{n} \tau \mu} (\mathbb{E}|W|^2)^{1/2} \left(\mathbb{E}\left|\frac{\overline{N_n}-n\nu}{\sqrt{\overline{N_n}}+\sqrt{n \nu}} \right|^2\right)^{1/2}
     &\leq& \frac{C\sigma''}{\sqrt{n} \tau \mu} \left(\mathbb{E}\left|\frac{N_n-n\nu}{\sqrt{n \nu}} \right|^2 \right)^{1/2} \nonumber\\
     &\leq& \frac{C\sigma''}{\sqrt{n\nu} \mu}.
 \end{eqnarray}

\noindent since $\mathbb{E}[W^2] =1 $ and $Var(N_n)=n \tau^2$. Also
note that for the second inequality we used
$|\overline{N}_n-n\nu|\leq |N_n- n\nu|$ which easily from the
definition of the truncation.

For the third term of the bound in (\ref{chenshaobound}), we have
\begin{eqnarray*}
  \sum_{i=1}^n \mathbb{E}|\xi_i(\Delta - \Delta_i)| \leq \sum_{i=1}^n (\mathbb{E}|\xi_i|^2)^{1/2} \mathbb{E} (|\Delta - \Delta_i|^2)^{1/2} &\leq& \sum_{i=1}^n \frac{1}{\sqrt{n}} \left(\mathbb{E}\left|\left(\frac{\sqrt{\overline{N_n}}-\sqrt{\overline{N_n-Y_i+\nu}}}{\sqrt{n}\tau \mu} \right)\sigma''Z_1\right|^2 \right)^{1/2}\\
&\leq& n \frac{\mathbb{E}|Z_1| \sigma''}{\sqrt{n}}
\left(\mathbb{E}\left|\frac{\sqrt{\overline{N_n}}-\sqrt{\overline{N_n-Y_1+\nu}}}{\sqrt{n}\tau
\mu} \right|^2 \right)^{1/2} \\
    &\leq& C \sqrt{n} \sigma'' \left(\mathbb{E} \left|\frac{\overline{N_n}-\overline{N_n-Y_1+\nu}}{\sqrt{n} \tau \mu (\sqrt{\overline{N_n}}+\sqrt{\overline{N_n-Y_1+\nu}})} \right|^2 \right)^{1/2}  \\
    &\leq& \frac{C \sigma''}{\tau \mu} \frac{(\mathbb{E}|Y_1-\nu|^2)^{1/2}}{\sqrt{n\nu/2}+ \sqrt{n \nu /2}}
\end{eqnarray*}

\noindent where we used $\mathbb{E}|\xi_i|^2 =1/n$, the identity
$\sqrt{x}-\sqrt{y} = \frac{x-y}{\sqrt{x}+\sqrt{y}}$ and the
inequality $|\overline{x}-\overline{x-y}|\leq |y|.$

We conclude
\begin{equation}\label{bounds3}
 \sum_{i=1}^n \mathbb{E}|\xi_i(\Delta - \Delta_i)|
=\sum_{i=1}^n \mathbb{E}\left|\frac{Y_i-\nu}{\sqrt{n \tau^2}}
(\Delta - \Delta_i) \right|  \leq \frac{C \sigma''}{\sqrt{\nu n}
\mu}.
\end{equation}

Using Theorem \ref{concentration}, we get

 \begin{eqnarray}
   \sup_{z \in \mathbb{R}} |\mathbb{P}(T_n'(Z_1) \leq z) - \mathbb{P}(T_n(Z_1,Z_2) \leq z) |  &\leq& \mathbb{P}(|N_n-n\nu|> n\nu/2) \nonumber\\
     &+&  \sup_{z \in \mathbb{R}}\mathbb{E}\big[|\mathbb{E}[\mathbbm{1}(\overline{T_n'}(Z_1) \leq z) - \mathbbm{1}(T_n(Z_1,Z_2) \leq z) ]\mathbbm{1}(|N_n-n\nu|\leq n\nu/2)| N_n\big] \nonumber\\
     &\leq& \sup_{z \in \mathbb{R}} |\mathbb{P}(W+ \Delta \leq z) -\mathbb{P}(Z_3 \leq
     z)| \nonumber\\
     &\leq& \frac{4\tau^2}{n \nu^2} + C \left(\frac{|Y_1|^3}{\tau^3\sqrt{n}}+ \frac{\sigma''}{\sqrt{n \nu} \mu}
     \right).
 \end{eqnarray}
where for the last step we combined the three estimates given in
(\ref{bounds1}), (\ref{bounds2}) and (\ref{bounds3}). Thus,

\begin{equation}\label{tn2}
    d_p(T_n'(Z_1),T_n(Z_1,Z_2)) \longrightarrow 0
\end{equation}
 as $n \rightarrow
\infty$ if  $\nu, \tau , \mu > 0$.

Finally we need to show that $d_P(T_n(Z_1), T_n'(Z_1))
\longrightarrow 0$. First observe that $$T_n(Z_1) - T_n'(Z_1) =
\frac{\sqrt{N_n} (\sigma' - \sigma'') Z_1}{\sqrt{n} b}
\longrightarrow 0$$ almost surely as $n \rightarrow \infty.$ Also we
know that $T_n'(Z_1)$ converges in distribution to $T_n(Z_1,Z_2)$.
Thus, using Slutsky's theorem we conclude that $T_n(Z_1)=T_n'(Z_1) +
T_n(Z_1)- T_n'(Z_1) $ also converges in distribution to
$T_n(Z_1,Z_2).$ Hence
\begin{equation}\label{tn3}
    d_P(T_n(Z_1),T_n'(Z_1))\leq d_P(T_n(Z_1),T_n(Z_1,Z_1))+
d_P(T_n(Z_1,Z_2), T_n'(Z_1)) \longrightarrow 0
\end{equation}
 as $n \rightarrow \infty.$

Hence combining (\ref{tn1}), (\ref{tn2}) and (\ref{tn3}), we obtain
$$d_p(T_n,T_n(Z_1,Z_2))\leq
d_p(T_n,T_n(Z_1))+d_p(T_n(Z_1),T_n'(Z_1))+d_p(T_n'(Z_1),T_n(Z_1,Z_2))\longrightarrow
0$$ as  $n \rightarrow \infty$ under the given assumptions and
result follows. \hfill $\square$

\vspace{0.15in}

\begin{remark}\label{prokhorovremark}
We can show that $d_P(T_n'(Z_1), T_n(Z_1,Z_2)) \rightarrow 0$ easily
if we are not interested in convergence rates. To see this, note
that we can write $T_n'(Z_1)$ as $$T_n'(Z_1) =
\sqrt{\frac{N_n}{n}}\frac{\sigma''}{b}\left(Z_1 + \frac{\frac{(N_n -
n\nu) }{\sqrt{n} \tau } \frac{\mu
\tau}{b}}{\sqrt{\frac{N_n}{n}}\frac{\sigma''}{b}} \right).$$ Now by
the strong law of large numbers $\frac{N_n}{n} \rightarrow \nu$ a.s.
and by the standard central limit theorem for independent random
variables $\frac{N_n - n\nu}{\sqrt{n} \tau} \rightarrow Z$ where $Z$
is a standard normal random variable independent of $Z_1$. Using
Slutsky's theorem twice with these observations immediately reveals
that $T_n'(Z_1)$ converges in distribution to a standard normal
random variable.
\end{remark}
\vspace{0.15in}

\textbf{Proof of Theorem \ref{chencase} : } First note that under
independence, we have $a_j = 0$ for $j=1,...,m$ so that $ \sigma' =
\sigma'' = \sigma.$ Following the proof of Theorem \ref{mainresult},
this implies that $d_K(T_n(Z_1), T_n'(Z_1)) = 0$  for every $n$. Now
the result follows from the estimates of $d_K(T_n,T_n(Z_1))$ and
$d_K(T_n(Z_1), T_n(Z_1,Z_2))$ by substituting $a_j = 0$ for
$j=1,...,m.$ \hfill $\square$

\vspace{0.15in}

\textbf{Proof of Corollary \ref{shang} :} In the proof of Theorem
\ref{mainresult}, we showed that $$d_K (H_n,Z_1) = d_K(T_n,
T_n(Z_1)) \rightarrow 0$$ where $H_n = \frac{\sum_{i=1}^{N_n}X_i -
N_n \mu}{\sqrt{N_n} \sigma'}$ and $(\sigma')^2 = \sigma^2 +2
\sum_{j=1}^m a_j \Gamma_{N_n,j}-\frac{2}{N_n} \sum_{j=1}^m j a_j.$
Since $\frac{\sigma'}{\sqrt{\sigma^2 + 2 \sum_{j=1}^m a_j}}
\longrightarrow 1 $ a.s., result follows from Slutsky's
theorem.\hfill $\square$

Finally we  give the proofs of the variance formulas given in
Proposition \ref{var1} and \ref{varrandom}.

\vspace{0.1in}

\textbf{Proof of Proposition \ref{var1} :} We have
\begin{eqnarray*}
                                                    Var \left(\sum_{i=1}^N X_i \right) &=& \sum_{i=1}^N Var(X_i) + 2 \sum_{1 \leq i <j \leq N} Cov(X_i, X_j) \\
                                                      &=& N \sigma^2+ 2 \sum_{j=1}^m (N-j) a_j \mathbbm{1}(N \geq j+1)
\end{eqnarray*}
Rearranging terms, we obtain
$$Var \left(\sum_{i=1}^N X_i \right) =  N \left(\sigma^2 + 2 \sum_{j=1}^m a_j \mathbbm{1}(N \geq j+1)  \right) - 2 \sum_{j=1}^m j a_j \mathbbm{1} (N \geq j+1)  $$

\noindent by which the variance formula follows. \hfill $\square$

\vspace{0.1in}

\textbf{Proof of Proposition \ref{varrandom} :} First note that
assumptions of Wald's identity are satisfied and so
$\mathbb{E}\left[\sum_{i=1}^{N_n} X_i \right] = n \nu \mu.$ Using
this, we get
\begin{eqnarray*}
  Var\left(\sum_{i=1}^{N_n} X_i \right) &=& \mathbb{E}\left[\left(\sum_{i=1}^{N_n} X_i  - n\nu \mu\right)^2 \right] \\
  &=&  \sum_{k=m+1}^{\infty} \mathbb{E}\left(\sum_{i=1}^k X_i -n \nu \mu \right)^2 \mathbb{P}(N_n = k) +  \sum_{k=0}^{m} \mathbb{E}\left(\sum_{i=1}^k X_i -n \nu \mu \right)^2 \mathbb{P}(N_n = k)
\end{eqnarray*}
where for the second equality we conditioned on $N_n$ which is
independent of $X_i'$s. Next note that we have
\begin{equation}\label{expec}
\mathbb{E} \left[\sum_{i=1}^k X_i \right] = k \mu \quad \text{and}
\quad \mathbb{E} \left(\sum_{i=1}^k X_i \right)^2 =k \sigma^2 + 2k
\sum_{j=1}^m a_j \Gamma_{j,k} - 2 \sum_{j=1}^m j a_j \Gamma_{j,k}
+k^2 \mu^2
\end{equation}

\noindent with $\Gamma_{j,k}=\mathbbm{1}(j \geq k+1)$. Thus, using
Proposition \ref{var1} and (\ref{expec}), and doing some elementary
manipulations, we obtain
\begin{eqnarray*}
  Var\left(\sum_{i=1}^{N_n} X_i \right) &=&    \sum_{k=m+1}^{\infty} \mathbb{E}\left(\sum_{i=1}^k X_i - k \mu + k \mu- n\nu \mu
   \right)^2 \mathbb{P}(N_n =k)\\
    &+& \sum_{k=0}^{m} \mathbb{E}\left[ \left(\sum_{i=1}^k X_i\right)^2 -2 n \nu \mu \left(\mathbb{E}\left[\sum_{i=1}^k X_i \right] \right) + n^2 \nu^2 \mu^2\right]\mathbb{P}(N_n
    =k)\\
    &=&  \sum_{k=m+1}^{\infty} \left(Var \left(\sum_{i=1}^k X_i \right) + (k \mu - n \nu
    \mu)^2
   \right)  \mathbb{P}(N_n =k) \\
   &+& \sum_{k=0}^{m} (k \sigma^2 + 2k \sum_{j=1}^m a_j \Gamma_{j,k} - 2 \sum_{j=1}^m j a_j \Gamma_{j,k} +k^2 \mu^2 -2 n \nu \mu^2 k + n^2 \nu^2 \mu^2
   ) \mathbb{P}(N_n =k)
\end{eqnarray*}

Noting that for $k \geq m+1$, $Var\left(\sum_{i=1}^k X_i\right) =
k\left(\sigma^2 +2 \sum_{j=1}^m a_j \right)- 2 \sum_{j=1}^m j a_j,$
we get

\begin{eqnarray*}
  Var\left(\sum_{i=1}^{N_n} X_i \right)   &=& \sum_{k=0}^{\infty} (k \sigma^2 + 2 k \sum_{j=1}^m a_j -2 \sum_{j=1}^m j a_j + k^2 \mu^2 -2 k n \nu \mu^2 + n^2 \nu^2
   \mu^2) \mathbb{P}(N_n =k) \\
   &+& \sum_{k=0}^m (k\sigma^2 + 2 k\sum_{j=1}^m a_j \Gamma_{j,k}- 2\sum_{j=1}^m j a_j \Gamma_{j,k} +\mu^2 k^2 -2 n \nu \mu^2 k + n^2 \nu^2 \mu^2
    \\
   &-& k \sigma^2 - 2k \sum_{j=1}^m a_j + 2 \sum_{j=1}^m j a_j - k^2 \mu^2 +2 k n \nu \nu^2 - n^2 \nu^2 \mu^2
   ) \mathbb{P}(N_n =k).
\end{eqnarray*}

After some cancelations and using the values for $\mathbb{E}[N_n]=n
\nu$ and $\mathbb{E}[N_n^2] = n \tau^2 + n^2 \nu^2$, we finally
arrive at

$$Var\left(\sum_{i=1}^{N_n} X_i \right) = n( \nu \sigma^2+ 2  \nu \sum_{j=1}^m a_j+ \mu^2 \tau^2) +
\alpha(m)$$ where $$\alpha(m) = \sum_{k=0}^m (2k \sum_{j=1}^m
a_j(\Gamma_{k,j}-1)-2 \sum_{j=1}^m j a_j (\Gamma_{k,j}-1)) -2
\sum_{j=1}^m j a_j.$$ The assertion that $\frac{\alpha(m)}{n}
\longrightarrow 0$  as $n \rightarrow \infty$ follows from the fact
that all the variables are bounded. \hfill $\square$

\section{Conclusion}

In this paper, we established a central limit theorem for random
sums of stationary $m-$dependent processes. Our proof is an
extension of the argument given in \cite{chenshaopaper} for the
i.i.d. case and this enables to recover their result. At the same
time, we were able to give variations of  the results in
\cite{shang}. In the subsequent research we are planning to (1)
obtain convergence rates for Theorem \ref{mainresult}, (2) relax the
$m-$dependence condition to a weak local dependence condition (For
such conditions, see \cite{chenshaolocal}), (3) adapt the size
biasing technique often used in normal approximation to the case of
random sums (See, for example, \cite{goldstein}) and (4) find more
applications on non-parametric statistics.

\end{document}